\documentclass{amsart}

\usepackage{enumerate}      
\usepackage{graphicx}       
\usepackage{url}            

\revertcopyright

\newtheorem{theorem}{Theorem}[section]
\newtheorem{lemma}[theorem]{Lemma}
\newtheorem{corollary}[theorem]{Corollary}
\newtheorem{algorithm}[theorem]{Algorithm}
\newtheorem{test}[theorem]{Test}

\theoremstyle{definition}
\newtheorem{definition}[theorem]{Definition}

\numberwithin{equation}{section}

\newcounter{c-save-section}
\newcounter{c-save-theorem}
\newcounter{c-sw-section}
\newcounter{c-sw-theorem}
\newcounter{c-nonhaken-section}
\newcounter{c-nonhaken-theorem}

\newcommand{\bdry}{\partial}
\newcommand{\co}{\colon\thinspace}
\newcommand{\inter}{\mathop{\mathrm{int}}}
\newcommand{\nbd}{\mathop{\mathrm{nbd}}}
\newcommand{\R}{\mathbb{R}}
\newcommand{\regina}{\emph{Regina}}
\newcommand{\swspace}{Weber-Seifert dodecahedral space}
\newcommand{\swurl}{\url{weber-seifert.rga}}
\newcommand{\tri}{\mathcal{T}}
\newcommand{\trisw}{\mathcal{T}_\mathrm{WS}}
\newcommand{\PS}{\mathcal{P}}
\DeclareMathOperator{\rank}{rk}
\DeclareMathOperator{\weight}{wt}

\begin{document}

%
\title[The Weber-Seifert dodecahedral space is non-Haken]
    {The Weber-Seifert dodecahedral space \\ is non-Haken}

\author{Benjamin A.~Burton}
\address{School of Mathematics and Physics \\
    The University of Queensland \\ Brisbane QLD 4072 \\ Australia}
\email{bab@maths.uq.edu.au}
\thanks{The first author is supported under the Australian Research Council's
Discovery funding scheme (project DP1094516).}

\author{J.~Hyam Rubinstein}
\address{Department of Mathematics and Statistics \\
    The University of Melbourne \\ VIC 3010 \\ Australia}
\email{rubin@ms.unimelb.edu.au}
\thanks{The second and third authors are partially supported under this same
scheme (projects DP0664276 and DP1095760).}

\author{Stephan Tillmann}
\address{School of Mathematics and Physics \\
    The University of Queensland \\ Brisbane QLD 4072 \\ Australia}
\email{tillmann@maths.uq.edu.au}


\begin{abstract}
In this paper we settle Thurston's old question of whether the
Weber-Seifert dodecahedral space is non-Haken, a problem that has been
a benchmark for progress in computational 3--manifold
topology over recent decades.
We resolve this question by combining recent
significant advances in normal surface enumeration, new heuristic
pruning techniques, and a new theoretical test that extends the seminal
work of Jaco and Oertel.
\end{abstract}

\maketitle

%
%

\section{Introduction} \label{s-intro}

\subsection{Motivation} \label{s-intro-mot}

In recent decades, 3--manifold topology has seen the resolution of many
important decision problems, from Haken's unknot recognition algorithm
in the early 1960s \cite{haken61-knot} through to the difficult and
multifaceted homeomorphism algorithm that was finally tied together
with Perelman's proof of the geometrisation conjecture in 2002.
Many of these resolutions, however, are in theory only.
Decision algorithms in 3--manifold topology are often highly complex and
extremely inefficient, and only a handful have ever
been implemented for practical use.  A key motivation in computational
topology is to elevate such algorithms from hypothetical procedures to
practical tools that can be employed in other theoretical applications.

In 1980, Thurston asked whether the Weber-Seifert dodecahedral space is
a Haken manifold \cite{birman80-problems}
(we define these concepts in Section~\ref{s-intro-defns}).
Resolving this question has been a symbolic goal for computational
topologists: in 1984, Jaco and Oertel devised an algorithm to decide
whether a 3--manifold is Haken \cite{jaco84-haken}, and the only barrier
to solving Thurston's question has been improving, implementing and running
the Jaco-Oertel algorithm.
Furthermore, resolving this specific question has broader implications---many
of the improvements to the Jaco-Oertel algorithm have a wider impact, since
this algorithm uses Haken's machinery of \emph{normal surface theory},
a toolset that has now become ubiquitous in 3--manifold decision problems.

Here we resolve the question of Thurston by proving the following theorem:

\setcounter{c-sw-section}{\arabic{section}}
\setcounter{c-sw-theorem}{\arabic{theorem}}
\begin{theorem}\label{t-sw}
The {\swspace} is non-Haken.
\end{theorem}

The proof is essentially computational; as noted earlier, the main
difficulty lies in redeveloping and extending the underlying algorithms
to the point where we can implement and execute them in a feasible timeframe.
We achieve this through a combination of techniques:
\begin{itemize}
    \item We implement and exploit recent advances in the core
    machinery of normal surfaces, including the Q-theory of
    Tollefson \cite{tollefson98-quadspace} and the enumeration
    algorithms of the first author \cite{burton09-convert,burton10-dd};

    \item We extend the work of Jaco and Oertel, essentially showing
    that we can reduce the computational workload by searching for
    compatible \emph{pairs} of normal surfaces instead of individual
    normal surfaces (Theorem~\ref{t-univ-disjoint} below);

    \item We develop new heuristic tests for identifying
    compressible surfaces; in theory these
    do not guarantee conclusive results, but in practice we find that
    they allow us to completely circumvent the most difficult step of
    the Jaco-Oertel algorithm (testing surfaces for compressibility).
\end{itemize}

\subsection{Definitions and key results} \label{s-intro-defns}

The closed, irreducible 3--manifold $M$ is \emph{Haken} if it contains
an embedded, injective surface different from the 2--sphere or
the projective plane. A surface $S$ in $M$ is \emph{injective} if the
inclusion $S \subset M$ induces a monomorphism between the fundamental
groups. The boundary of a regular neighbourhood of an injective surface
is 2--sided and (geometrically) incompressible.

The \emph{Weber-Seifert dodecahedral space} is formed by identifying
opposite faces of a dodecahedron with a $3/10$ twist, and it was one of
the first known examples of a hyperbolic 3--manifold
\cite{weber33-dodecahedral}. For our computations we use
a triangulation of this manifold with 23
tetrahedra. There are at least three distinct
triangulations with this number of tetrahedra, and they are the smallest
triangulations of the Weber-Seifert dodecahedral space known to the
authors.

The setting for this paper is the following. Let $M$ be a closed,
orientable, irreducible 3--manifold with a fixed triangulation $\tri$.
Work of Jaco and Oertel \cite{jaco84-haken}, Tollefson
\cite{tollefson95-isotopy} and Oertel \cite{oertel86-branched} studies
the set of isotopy and projective isotopy classes of closed, injective
surfaces in $M$ and how they are represented in the projective solution
space $\PS$ of normal surface theory. Here, two embedded,
2--sided surfaces in $M$ are in the same \emph{projective isotopy class}
if there exist multiples of each which are isotopic.

A natural environment for algorithmic topology on a 3--manifold is a
\emph{0--efficient triangulation} \cite{jaco03-0-efficiency}.
A triangulation of a closed 3--manifold is 0--efficient if the only
normal 2--spheres are the links of the vertices of the triangulation
(that is, frontiers of small regular neighbourhoods of these vertices).
Standard facts about normal surfaces are recalled in
Section~\ref{s-universal}, and the only non-standard term that needs to
be clarified is the following: A surface in $M$ is termed a
\emph{vertex surface} if it is a connected, 2--sided normal
surface and the ray from the origin through its normal coordinate
passes through a vertex of $\PS$. Both the nomenclature and the
definition of this concept vary widely in the literature; the present is
chosen so that every vertex surface is 2--sided and there
is a unique vertex surface associated with each
admissible vertex of $\PS$ (an \emph{admissible} vertex has at most one
non-zero quadrilateral coordinate per tetrahedron).
In this terminology, Jaco and Oertel
\cite{jaco84-haken} show that $M$ contains an embedded, injective surface
if and only if one of the vertex surfaces is incompressible.
(As usual, a sphere is not incompressible.)

In Section~\ref{s-universal} we prove the following key result regarding
vertex surfaces:

\begin{theorem} \label{t-univ-disjoint}
Suppose that the closed, orientable, irreducible 3--manifold $M$ is Haken,
and let $\tri$ be a 0--efficient triangulation of $M$.
Then one of the cases 1~and~2 holds:
\begin{enumerate}[1.]
    \item \emph{(General case)}
    There are two
    distinct, incompressible vertex surfaces which are compatible.
    \item \emph{(Exceptional case)}
    There are at most finitely many projective isotopy classes of
    injective surfaces in $M$, and at least one of the following holds:
    \begin{enumerate}[(a)]
        \item the manifold $M$ fibres over the circle and there is a
        unique projective isotopy class of injective surfaces---in
        particular, $\rank H^1(M; \R) =1$ and the class is represented
        by a fibre; or
        \item there is a solid torus in $M$ which contains at least
        three edges of $\tri$ and has normal boundary.
    \end{enumerate}
\end{enumerate}
\end{theorem}

Note that, even if there are at most finitely many projective isotopy classes
of injective surfaces in $M$, the triangulation may also satisfy the general
case.

The following
corollary follows immediately from the above result. It is also a
consequence of work of Oertel \cite{oertel86-branched} and Tollefson
\cite{tollefson95-isotopy}, where faces of the projective solution
space are studied using branched surfaces.

\begin{corollary}
Let $M$ be a closed, orientable, irreducible 3--manifold, and let $\tri$ be
a 0--efficient triangulation of $M$. If $\rank H^1(M; \R) \ge 2$, then
there are two distinct, incompressible vertex surfaces which are compatible.
\end{corollary}

In Section~\ref{s-heuristic} we introduce some simple tests that can
help identify when a 2--sided surface is compressible. These tests are
merely heuristic techniques---there is no guarantee for any particular
surface that they will give a conclusive result, nor can they prove a
surface to be incompressible. However, these techniques
are found to be surprisingly
effective in practice.  Indeed, combined with some human intervention and the
original algorithm of Jaco and Oertel, they can completely resolve the
question of whether the {\swspace} is non-Haken (Theorem~\ref{t-sw}),
as noted at the end of Section~\ref{s-swproof}.

Combining Theorem~\ref{t-univ-disjoint} with these heuristic techniques,
we obtain the following new test to identify whether the given manifold $M$
is non-Haken:

\setcounter{c-nonhaken-section}{\arabic{section}}
\setcounter{c-nonhaken-theorem}{\arabic{theorem}}
\begin{test}[Is $M$ non-Haken?] \label{test-nonhaken}
First check that $\rank H^1(M; \R) =0$, since $M$ is otherwise Haken. Then
compute a 0--efficient triangulation and enumerate all vertex surfaces.
Check that each vertex surface of zero Euler characteristic misses at most
two edges. If this is not the case, the test is inconclusive. Otherwise,
let $\mathcal{S}=\emptyset$. For each vertex surface, either determine a
compression disc (for instance, using the heuristic techniques mentioned
above), or add it to the set $\mathcal{S}$. Last, check that no
two members of $\mathcal{S}$ are compatible. If this is the case, $M$ is
non-Haken; otherwise the test is inconclusive.
\end{test}

This na\"ive test suffices to prove Theorem~\ref{t-sw}.  The
projective solution space of the chosen triangulation with 23 tetrahedra has
1751 admissible vertices: one is represented by a 2--sphere which links
the vertex,
24 are represented by tori which link the edges, and the remaining 1726
vertices are represented by surfaces of negative Euler characteristic.
Heuristic pruning finds compressing discs for \emph{all but 16} of them,
and no two of the remaining surfaces are compatible. The total computation
time was a little over six hours; details are given in Section~\ref{s-swproof}.

Other applications of Theorem~\ref{t-univ-disjoint}
are given in Section~\ref{s-universal}. All routines
developed in this paper are implemented in version~4.6.1 of the
open-source software package {\regina} \cite{regina,burton04-regina}.

\section{Compatible, injective vertex surfaces} \label{s-universal}

The original algorithm of Jaco and Oertel \cite{jaco84-haken} to
determine whether a manifold is Haken requires the enumeration of all
vertex surfaces, and for each such surface one needs to check whether it
is incompressible or not. Deciding incompressibility is computationally
very expensive, and so this section uses the work of Jaco and Oertel to improve
on their algorithm.


\subsection{Triangulation}

The notation and terminology of \cite{jaco03-0-efficiency} will be used
in this paper. Hence a triangulation $\tri$ consists of a union of $t$
pairwise disjoint 3--simplices, $\widetilde{\Delta}$, a set of face
pairings, $\Phi$, and a natural quotient map $p\co \widetilde{\Delta}
\to \widetilde{\Delta} / \Phi = M$. This is often referred to as a
\emph{semi-simplicial} or \emph{singular} triangulation since not all
simplices are necessarily embedded in $M$. The space
$\widetilde{\Delta}$ has a natural simplicial structure with four
vertices for each 3--simplex. The quotient map $p$ is required to be
injective on the interior of each simplex of each dimension. The image
of a simplex in $\widetilde{\Delta}$ under $p$ is a singular simplex in
$M$. It is customary to refer to the image of a 3--simplex as a
\emph{tetrahedron in $M$} (or \emph{of the triangulation}) and to refer
to its faces, edges and vertices with respect to the pre-image.
Similarly for images of 2--, 1--, and 0--simplices, which will be
referred to as \emph{faces}, \emph{edges} and \emph{vertices in $M$} (or
\emph{of the triangulation}) respectively. If a singular simplex is
contained in $\partial M$, then it is termed \emph{boundary} (such as
a \emph{boundary edge} or a \emph{boundary face}); otherwise
it is termed \emph{internal}. Notice that an internal singular simplex
need not be disjoint from $\partial M$. A \emph{normal isotopy of $M$}
is an isotopy which leaves the image of the interior of every simplex in
$\widetilde{\Delta}$ invariant. The quotient space $M$ is a manifold if
the link of each vertex in $M$ is a sphere or a disc.


\subsection{Normal surfaces}

This terminology again follows \cite{jaco03-0-efficiency}. A normal
surface in the triangulated 3--manifold $M$ meets every tetrahedron in a
pairwise disjoint, finite
union of discs which are \emph{normal triangles} or \emph{normal
quadrilaterals}. A normal surface is hence a properly embedded surface
in $M$. The \emph{normal coordinate} is a point in $\R^{7t}$ that
records the number of discs of each type
in a normal surface. It satisfies a system of integral, linear equations,
termed the \emph{matching equations}. The set of all solutions with
non-negative coordinates to this system is intersected with the affine
subspace consisting of all points whose coordinates sum to one to give
the \emph{projective solution space} $\PS$. This is a (bounded) polytope
whose vertices have rational coordinates.
Given any normal surface, its normal coordinate determines a
unique point in $\PS$.

A point in $\R^{7t}$ is \emph{admissible} if all of its coordinates are
non-negative and at most one quadrilateral coordinate from each
tetrahedron is non-zero. Each integral admissible solution to the
matching equations determines a unique normal surface and vice versa.
Two normal surfaces are said to be \emph{compatible} if they do not meet
a tetrahedron in quadrilateral discs of different types. This is the
case if and only if the sum of their normal coordinates is admissible.

A surface in $M$ is termed a \emph{vertex surface} if it is a
connected, 2--sided normal surface and the ray from the origin through
its normal coordinate passes through a vertex of $\PS$. Any ray from
the origin through an admissible vertex of $\PS$ contains a unique vertex
surface.
A vertex surface is sometimes termed a \emph{vertex normal surface} or a
\emph{fundamental edge surface} in the literature. It should not be
confused with a \emph{vertex linking surface} (which is the boundary of
a small neighbourhood of a vertex in the triangulation).


\subsection{$0$--efficiency}

The triangulation $\tri$ is \emph{$0$--efficient} if the only normal
2--spheres are vertex linking. It is shown in \cite{jaco03-0-efficiency}
that any triangulation of a closed, orientable, irreducible 3--manifold
can be modified to a $0$--efficient triangulation unless $M = S^3, \R P^3$
or $L(3,1)$. Moreover, the conversion algorithm is implemented in
{\regina} \cite{burton04-regina}. The algorithm typically takes only
marginally longer than the time required to
enumerate the admissible vertex
solutions to the so-called $Q$--matching equations of Tollefson
\cite{tollefson98-quadspace}. Detailed time trials can be found in
\cite{burton10-dd}.


\subsection{The work of Jaco and Oertel}

Jaco and Oertel \cite{jaco84-haken} give an algorithm to decide whether
a triangulated manifold is Haken. Surfaces are analysed using handle
decompositions, but since their arguments are topological, the results
carry over to triangulations. In this subsection, one of their key
results will be re-stated in a topological version. This version follows
verbatim from the proof of Theorem~2.2 in \cite{jaco84-haken}.

Let $M$ be a closed, irreducible 3--manifold, and $F_1$ and $F_2$ be
embedded surfaces in general position. Then $F_1 \cap F_2$ is a
finite union of pairwise disjoint curves. A component  of $F_1 \cap F_2$
is termed a \emph{switch curve}. Let $\gamma$ be a switch curve. A
regular neighbourhood $N(\gamma)$ of $\gamma$ is chosen such that
$N(\gamma) \cap F_i$ is a regular neighbourhood of $\gamma$ in $F_i$ for
each $i$. Since $N(\gamma)$ is either a solid Klein bottle or a solid
torus, it follows that $\partial N(\gamma) \setminus (F_1 \cup F_2)$
consists either of two M\"obius bands or of four annuli, called
\emph{switch bands} or \emph{switch annuli} respectively. Two switch
annuli are said to be \emph{opposite} if they do not share a boundary
component. A \emph{switch along $\gamma$} consists of deleting the
portion of $F_1 \cup F_2$ inside $N(\gamma)$ and connecting the free
boundary components by a switch band or by two opposite switch annuli
(depending on whether $\partial N(\gamma)$ is a Klein bottle or torus
respectively). It follows that there are two possible switches along
$\gamma$.

Denote by $F_1 + F_2$ the surface obtained from $F_1 \cup F_2$ where at
each component of $F_1 \cap F_2$ one of the two possible switches has
been chosen. If $F_1$ and $F_2$ are compatible normal surfaces with
respect to a handle decomposition or a triangulation, then there is a natural
choice at each switch curve, called a \emph{regular switch}, such
that $F_1 + F_2$ is again a normal surface. 

The surface $F=F_1 + F_2$ is said to be in \emph{reduced form} if it cannot
be written as $F=F'_1 + F'_2$, where $F'_i$ is isotopic to $F_i$ in $M$
and $F'_1 \cap F'_2$ has fewer components than $F_1 \cap F_2$. It should
be noted that in these two sums, the embedding of $F$ in $M$ is the same
(these are not equalities up to isotopy), and that any sum can be
changed to a sum in reduced form.

We will also denote by $F_1 +_\gamma F_2$ the surface obtained from
$F_1 \cup F_2$ by choosing the same switches as for $F_1 + F_2$ except for
the curve $\gamma$, where the other switch possibility is chosen. 

\begin{theorem}[Jaco-Oertel]\label{t-top-jo}
Let $M$ be a closed, irreducible 3--manifold and $F$ be an embedded, 2--sided
and incompressible surface. If $F=F_1 + F_2$ is in reduced form, then either
\begin{enumerate}
\item $F_1$ and $F_2$ are incompressible, or
\item there exists $\gamma \in F_1 \cap F_2$ such that the surface
$F_1 +_\gamma F_2$ has two components, $F'$ and $T$, with the property
that $T$ is a torus which bounds a solid torus in $M$ and has longitude
isotopic to $\gamma$, and $F'$ is isotopic to $F$ via an isotopy which
fixes the complement of a neighbourhood of the union of the solid torus
and $N(\gamma)$. Moreover, $T = A \cup A'$, where $A'$ is a switch annulus
not contained in $F$ and $A$ is an annulus contained in $F$.
\end{enumerate}
\end{theorem}

The second possibility is illustrated in Figures 6 and 8 of
\cite{jaco84-haken}. Under the additional assumption that $F$ is a
normal surface of least weight (with respect to a triangulation or
handle decomposition), one sees that the second alternative is not
possible and concludes that the summands are incompressible. Under these
circumstances, one can also omit the hypothesis that the surface be in
reduced form. This is the result stated in \cite{jaco84-haken}.


\subsection{Proof of Theorem \ref{t-univ-disjoint}}

Since $M$ is Haken, the $0$--efficient triangulation $\tri$ has a single
vertex, $v$. Suppose $S$ is a connected, injective surface in $M$. By
possibly replacing $S$ with the boundary of a regular neighbourhood of
$S$, one may assume that $S$ is a 2--sided, geometrically incompressible
surface. Recall that the weight of an embedded surface in $M$
which is in general position with respect to the triangulation is the
cardinality of its intersection with the 1--skeleton, written
$\weight (F) = | \ F \cap \tri^{(1)} \ |$. After performing an isotopy, one
may also assume that $S$ is a normal surface of least weight in its isotopy
class. If $S$ is not a vertex surface then, as in
\cite[Corollary~3.4]{jaco84-haken}, it follows from Theorem \ref{t-top-jo}
that there are two compatible, injective vertex solutions. Since $\PS$ has
finitely many vertices, it follows that if there are infinitely many
projective isotopy classes of injective surfaces in $M$, then there are
two compatible, injective vertex solutions. This is the
general case (Case~1) of the theorem and it remains to show that otherwise
we are in the exceptional case (Case~2a or~2b).

Hence suppose that there are at most finitely many projective isotopy
classes of injective surfaces in $M$ and that $S$ is a vertex surface.
We first produce a second normal surface isotopic but not normally
isotopic to $S$ by performing a \emph{finger move}. That is, we push a
portion of the surface along an edge and across the vertex in a
controlled fashion.  The details are as follows.

Let $e$ be an edge of the triangulation. An intersection point
$p \in S \cap e$ is \emph{outermost} if one of the two components
of $e \setminus \{  p, v\}$ does not contain any other intersection points.
Since $S$ is not a vertex linking sphere, it follows that there is an
edge $e_0$ having an outermost point $p_0 \in S \cap e_0$ which is
incident with a quadrilateral disc in $S$. Orient $e_0$ such that
travelling from $p_0$ in the positive direction to $v$, one does not meet $S$.
Since $S$ is 2--sided, give $S$ a transverse orientation which agrees at
$p_0$ with the orientation of $e_0$. Let $N(e_0)$ be a small regular
neighbourhood of $e_0$. Since $M$ is orientable, this is a torus,
and $S \cap N(e_0)$ consists of meridian discs since $S$ is normal.
Denote by $D$ the connected component of $S \cap N(e_0)$ passing through
$p_0$. We now perform an isotopy of $S$ which fixes $S \setminus D$ and
moves $p_0$ along $e_0$ in the positive direction just past $v$.
The resulting surface, $S_1$, will not be normal and can be chosen such that
\[\weight (S_1) = \weight (S) + 2(E-1),\]
where $E$ is the number of edges of the triangulation. Now push $S_1$
slightly off $S$ in the positive direction wherever they agree, giving a
surface $S_2$ which is disjoint from $S$ and has the same weight as $S_1$.

Since no face in the triangulation is a cone \cite{jaco03-0-efficiency}
or a dunce hat \cite{jaco09-minimal-lens}, and $S$ is pushed in the direction
of the oriented edge $e_0$, it follows from analysing the resulting isotopy of a
quadrilateral disc in $S$ meeting $p_0$, that there is at least one face
of the triangulation which $S_2$ meets in a return arc. Since $S_2$ is
incompressible, it can be normalised by isotopies giving a normal surface
$S_3$ which satisfies:
\[\weight (S_3) \le \weight (S) + 2(E-2).\]
Note that $S_3$ is disjoint from $S$, since $S$ acts as a barrier surface
for the normalisation of $S_2$. If $S_3$ is normally isotopic to $S$, then
there is a product region between the two surfaces which does not
contain the unique vertex of the triangulation. However, there also is a
product region between $S_3$ and $S$ containing the vertex which arises
from the first isotopy. It follows that $M$ fibres over the circle with
fibre $S$. If every connected, injective surface in $M$ is the fibre in
some fibration over the circle then there is a unique projective isotopy
class of injective surfaces in $M$ (since we assume that there are at
most finitely many). This is the first exceptional case stated in the theorem,
Case~2a.

Hence assume that $S$ is not a fibre in a fibration of $M$ over the
circle, which means that $S_3$ is not normally isotopic to $S$. If $S_3$ is a
vertex surface we have the general case of the theorem.
Note also that, by our assumption that $S$ is least weight in its isotopy
class, it follows that $\weight (S_3) \ge \weight (S)$.
If $\weight (S_3) = \weight (S)$, then we can apply  \cite{jaco84-haken} to
conclude that either $S_3$ is a vertex surface or
$S_3$ is a sum of vertex surfaces which are all incompressible.
In the latter case, each of the vertex surfaces is compatible with $S$
and they cannot all be copies of $S$ and
so we again have the general case of the theorem.

Hence suppose that $\weight (S_3) > \weight (S)$, and that $S_3$ is not a
vertex surface. We may further assume that $S_3$ has minimal weight amongst
all normal surfaces that are \emph{disjoint} from $S$ and are isotopic to $S$
but not normally isotopic to $S$. Since $S_3$ is not a vertex surface,
there is a positive integer $n$ and there are vertex surfaces $V_i$ and
positive integers $n_i,$ such that 
\[n S_3 = \sum n_i V_i,\]
where $m F$ signifies $m$ pairwise disjoint, parallel copies of $F$
(each normally isotopic to $F$), and the sum uses the usual
\emph{regular switches} with respect to the triangulation. Since $S_3$ is
compatible with $S$, each $V_i$ is compatible with $S$. Moreover, at least
one $V_i$ is not normally isotopic to $S$ since otherwise $n S_3$ is isotopic
to $(\sum n_i) S$. We can therefore write 
\[n S_3 = V + W,\]
where $V$ is a vertex surface distinct from but compatible with $S$.
There are normal surfaces $V'$ and $W'$ (isotopic in $M$ to $V$ and $W$
respectively), such that $n S_3=V' + W'$ in reduced form. 

According to Theorem \ref{t-top-jo}, we have two cases to consider. In
the first case, $V'$ (and hence $V$) is incompressible, and therefore
$S$ and $V$ are two distinct compatible, injective, vertex surfaces.
This is the general case of the theorem.

In the second case, following  \cite{jaco84-haken},
there exists $\gamma \in V' \cap W'$ such that the
surface $V' +_\gamma W'$ is the disjoint union of two surfaces, one of
which is isotopic to $n S_3$ and the other is a compressible torus $T$. Since the
switch curve $\gamma$ is 2--sided in $V'$ and $W'$, it follows that
$V'+_\gamma W'$ has $(n-1)$ components which are normally isotopic copies
of $S_3$. In addition, there is the torus $T$ and a component
$X$ which is isotopic to $S_3$. The latter normalises to give a surface $S_4$.

If $S_4$ is disjoint from $S$, then it must be normally isotopic to $S$
since its weight is strictly less than the weight of $S_3$. This gives an
isotopy from $S_3$ to $S$ which must pass through the vertex since otherwise
$M$ is a bundle. It follows that the vertex is contained in the solid torus
bounded by $T$. Considering weight, we have:
\[n \weight (S_3) > (n-1) \weight (S_3) +  \weight (S) + \weight (T),\]
whence 
\[\weight (S) + 2(E-2) \ge \weight (S_3) > \weight (S) + \weight (T),\]
giving 
\[2(E-2) > \weight (T).\]
Since $T$ is separating, it must meet each edge an even number of times,
whence the solid torus bounded by it contains at least three edges. The
boundary of a regular neighbourhood of one of the edges is an embedded torus
and a barrier surface (see \cite{jaco03-0-efficiency}).
Since $M$ is irreducible, the process of normalising $T$ 
in the complement of this torus shows that $T$ either
shrinks to a normal surface using isotopies or it shrinks to a 2--sphere embedded 
in a tetrahedron using isotopies and a single compression.
In the latter case, $M$ is a lens space, contradicting the fact that it is
Haken. Since normalisation does not increase the weight, we have the second
exceptional case, Case~2b, stated in the theorem.

It remains to analyse the possibility that $S_4$ is not disjoint from
$S$. In this case, $X$ cannot be disjoint from $S$ since otherwise $S$
is a barrier for the normalisation of $X$ to $S_4$. But if $S \cap X
\neq \emptyset$, then $S$ meets a neighbourhood of the union of the
solid torus bounded by $T$ and $N(\gamma)$ in a union of annuli and
can be isotoped to be disjoint from this and $X$.
(See Figure 8 of \cite{jaco84-haken} for an illustration
of the isotopy.) This would reduce the weight of $S$, which is a
contradiction.
Therefore $S$ does not meet $X$ and hence does not meet $S_4$. This
concludes the proof of the theorem. \qed


\subsection{Applications of Theorem \ref{t-univ-disjoint}}

Throughout this subsection, $M$ denotes a closed, orientable,
irreducible 3--manifold with fixed $0$--efficient triangulation $\tri$.

\begin{definition}[Large normal torus]
A normal torus disjoint from at least three edges will be termed a
\emph{large normal torus}.
\end{definition}

For instance, the boundary of a regular neighbourhood of a layered solid
torus subcomplex in $M$ shrinks to a large normal torus unless $M$ is a
lens space (see, for instance, \cite{burton04-facegraphs,jaco09-minimal-lens}
for a definition of this subcomplex).
Such a subcomplex appears in the natural triangulations of Dehn
fillings of knot complements. Large normal tori therefore occur in many natural
triangulations of both Haken and non-Haken manifolds.

\begin{algorithm}[Large normal torus recognition]\label{a-largetorus}
To check whether there is a large normal torus in $M$, it is necessary and
sufficient to verify that each vertex surface of Euler characteristic zero
is disjoint from at most two edges.
\end{algorithm}

This algorithm follows immediately from the fact that (i)~edge
weights and Euler characteristic are additive, and (ii)~the only normal
2--spheres
are vertex linking and hence can be made disjoint from any other normal
surface.

If $M$ is atoroidal and $\tri$ is $0$--efficient, then every normal torus
bounds a solid torus \cite{jaco09-minimal-lens}. It follows that in this case,
Algorithm~\ref{a-largetorus} can be used
to decide whether there is a solid torus in $M$
which contains at least three edges and has normal boundary.
In the general case, Algorithm~\ref{a-largetorus}
only helps in certain cases to identify when such a solid torus does not
exist.

This is also the philosophy in the tests below. Whilst in theory, there are
algorithms to determine whether a given manifold is non-Haken or a tiny
bundle, they often turn out to be impractically slow. Below are some simple
tests that allow an answer in a feasible amount of time, even though
they may not always be conclusive. We begin with the new non-Haken
test described in the introduction.

\setcounter{c-save-section}{\arabic{section}}
\setcounter{c-save-theorem}{\arabic{theorem}}
\setcounter{section}{\arabic{c-nonhaken-section}}
\setcounter{theorem}{\arabic{c-nonhaken-theorem}}
\begin{test}[Is $M$ non-Haken?]
First check that $\rank H^1(M; \R) =0$, since $M$ is otherwise Haken. Then
enumerate all vertex surfaces. Use Algorithm~\ref{a-largetorus} to
check that there is no large normal torus. If there is a large normal torus,
the test is inconclusive. Otherwise, let $\mathcal{S}=\emptyset$. For each
vertex surface, either determine a compression disc (for instance, using
the heuristic techniques described in Section~\ref{s-heuristic}),
or add it to the set
$\mathcal{S}$. Last, check that no two members of $\mathcal{S}$ are compatible.
If this is the case, then $M$ is non-Haken; otherwise the test is inconclusive.
\end{test}
\setcounter{section}{\arabic{c-save-section}}
\setcounter{theorem}{\arabic{c-save-theorem}}

Several routines are described in Section~\ref{s-heuristic} that search for
compression discs that are computationally easy to find.  These routines are
often sufficient to keep the list $\mathcal{S}$ relatively short. There
are, however, many classes of triangulated non-Haken manifolds which the
above approach will not recognise as non-Haken. For instance, if the
triangulation contains a layered solid torus subcomplex, then there is a large
normal torus as noted earlier.

\begin{definition}[Tiny bundle]
The closed, orientable, irreducible 3--manifold $M$ is a \emph{tiny bundle}
if $M$ fibres over the circle and there is a unique projective isotopy class
of injective surfaces in $M$.
\end{definition}

\begin{test}[Is $M$ a tiny bundle?]
The computational steps are as in Test \ref{test-nonhaken}, except that one
first checks that $\rank H^1(M; \R) =1$.
\end{test}

\begin{test}[Is there a finite number of projective isotopy classes?]
First show that $M$ is Haken. Then enumerate all vertex surfaces. Let
$\mathcal{S}=\emptyset$. For each vertex surface, either determine a
compression disc, or add it to the set $\mathcal{S}$. Last, check that no
two members of $\mathcal{S}$ are compatible. If this is the case, then
there are at most finitely many projective isotopy classes of injective
surfaces; otherwise the test is inconclusive.
\end{test}

\section{Heuristic pruning} \label{s-heuristic}

In this section we introduce some simple tests that can help
identify when an embedded surface within a closed 3--manifold is
compressible.  These are merely heuristic techniques---there
is no guarantee for any particular surface that they will give a conclusive
result.  Moreover, they work in one direction only---they can never show
a surface to be \emph{in}compressible.

With these tests, we are able to take a list of \emph{potential}
incompressible surfaces (such as the vertex normal surfaces within a
triangulation) and filter out irrelevant surfaces from this list.  This
leaves us fewer surfaces on which we must run more expensive
procedures, such as the conclusive but extremely slow
incompressibility algorithm of Jaco and Oertel \cite{jaco84-haken}.

The key idea behind these heuristic tests is to search within a bounded
3--manifold triangulation for embedded discs with simple
combinatorial structures.  Lemmata~\ref{l-prune-face}
and~\ref{l-prune-equator} describe these structures, and
Algorithm~\ref{a-heuristic} shows how they can be used effectively to
test for compressibility.  Although these tests are simple in theory,
we see in Section~\ref{s-swproof} of this paper that they can be
surprisingly effective in practice.

\begin{lemma} \label{l-prune-face}
    Let $\tri$ be a triangulation of a bounded 3--manifold $M$
    (that is, a compact 3-manifold with non-empty boundary).
    Let $F$ be a non-boundary face of $\tri$, and suppose that
    all three edges of $F$ lie entirely within the boundary $\bdry M$,
    as illustrated in the leftmost diagram of Figure~\ref{fig-prune-face}.
    Note that neither the edges of $F$ nor the vertices of $F$ are
    required to be distinct.

    \begin{figure}
    \centering
    \includegraphics{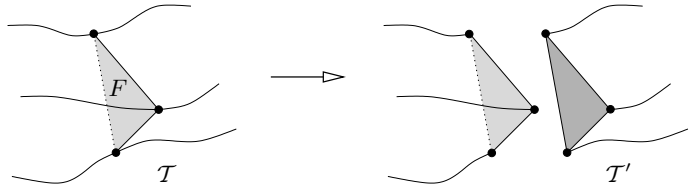}
    \caption{Cutting along a face with three boundary edges}
    \label{fig-prune-face}
    \end{figure}

    Suppose that we ``unglue'' the two tetrahedra on either side of $F$
    (that is, we remove the corresponding pair of tetrahedron faces
    from the list of face identifications that make up $\tri$), as
    illustrated in the rightmost diagram of Figure~\ref{fig-prune-face}.
    Then the result is a new triangulation $\tri'$ of some 3--manifold $M'$,
    which is homeomorphic to $M$ sliced along a properly embedded disc.
\end{lemma}

Note that we do not describe this operation as ``slicing $\tri$ along the
face $F$'', since $F$ might have self-intersections and therefore
might not be embedded.  However, whether or not $F$ is embedded,
the act of ungluing the two tetrahedra on either side of $F$ is well-defined
and simple to perform.  Proving that self-intersections of $F$ do not matter
is in fact the main point of this lemma.

\begin{proof}
    If the edges and vertices of $F$ are all distinct, then this result
    is straight\-for\-ward---the face $F$ forms a properly embedded disc in
    $M$, and the new triangulation $\tri'$ is just $M$ sliced along this
    disc.

    Consider then the case where different edges and/or vertices of $F$ are
    identified.  We first prove that $\tri'$ is indeed a 3--manifold
    triangulation, and then we show that the corresponding manifold $M'$ has
    the required property.

    The only situations in which $\tri'$ might \emph{not} be a 3--manifold
    triangulation are (i)~where some edge of $\tri'$ is identified with itself
    in reverse, and (ii)~where some vertex of $\tri'$ does not have a small
    closed neighbourhood that is a 3--ball.

    The first situation is easily eliminated, since ungluing the tetrahedra
    on either side of $F$ cannot create any new edge identifications (though
    it can remove them).  Consider then some vertex $V$ of the face $F$
    in $\tri$, and let $\nbd(V)$ be a small closed neighbourhood of $V$.
    Since the edges and vertices of $F$ all lie on the boundary $\bdry M$,
    the neighbourhood
    $\nbd(V)$ must be a 3--ball with $V$ on its boundary, as illustrated in
    the leftmost diagram of Figure~\ref{fig-vertex-slice}.

    \begin{figure}
    \centering
    \includegraphics{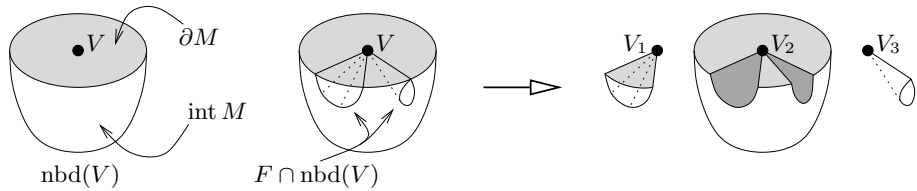}
    \caption{The ungluing operation in the vicinity of a vertex of $F$}
    \label{fig-vertex-slice}
    \end{figure}

    Now consider what happens when we unglue the tetrahedra on either side
    of $F$.  The intersection $F \cap \nbd(V)$ consists of up to three
    ``triangular'' discs in $\nbd(V)$ (one for each corner of $F$ that
    meets $V$).  Note that these discs might be joined along sections of
    their boundaries (corresponding to edges of $F$) to form larger
    discs or even branched structures.
    Examples of such discs are shown in the centre diagram of
    Figure~\ref{fig-vertex-slice}.

    The interiors of these individual triangular discs are
    disjoint and embedded in $\inter(\nbd(V))$, and the
    boundaries of these discs lie in the boundary $\bdry \nbd(V)$.
    Although the disc boundaries might intersect on $\bdry M$ (as a result of
    identifications between edges of $F$), they can never intersect in
    $\inter(M)$---in other words, the portions of the disc
    boundaries within $\inter(M)$ are also disjoint and embedded.

    We find that, when we unglue the tetrahedra on either side of $F$,
    we effectively slice $\nbd(V)$ along these discs, as illustrated in the
    rightmost diagram of Figure~\ref{fig-vertex-slice}.  This divides
    $\nbd(V)$ into several smaller 3--balls, splitting $V$ into several
    different vertices as a consequence.  Note that any disc edges
    that are pinched together on $\bdry M$ will fall apart, since
    only the face $F$ was holding them together (this happens
    with the rightmost disc in Figure~\ref{fig-vertex-slice}).

    Because of the well-behaved manner in which these discs are placed within
    $\nbd(V)$, we see that every resulting vertex of $\tri'$ has a small
    closed neighbourhood that is a 3--ball (that is, no ``bad'' holes have
    been cut out of $\nbd(V)$ as illustrated in Figure~\ref{fig-bad-hole}).
    Therefore $\tri'$ is indeed a 3--manifold triangulation, and we denote
    the corresponding 3--manifold by $M'$.
    Figure~\ref{fig-prune-pinched} illustrates the entire transformation
    from $\tri$ to $\tri'$ in the case where all three vertices of $F$
    are identified (here the single vertex $V$ in $\tri$ splits into four
    vertices in $\tri'$).

    \begin{figure}
    \centering
    \includegraphics{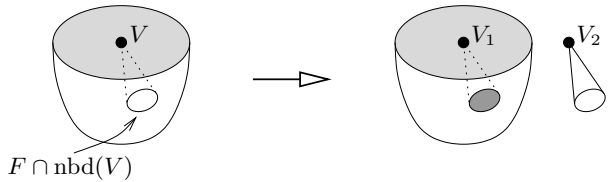}
    \caption{A ``bad'' hole where $V$ has no 3--ball neighbourhood}
    \label{fig-bad-hole}
    \end{figure}

    Now that we know that $\tri'$ is indeed a triangulation of the 3--manifold
    $M'$, the remainder of the lemma is straightforward to prove.
    Define $N$ to be the 3--manifold obtained by removing a small open
    neighbourhood of the boundary $\bdry M$ from $M$.  Likewise,
    let $N'$ be obtained by removing a small open neighbourhood of $\bdry M'$
    from $M'$.  It is clear that $F \cap N$ is a properly embedded disc
    in $N$, and that the 3--manifold $N'$ is obtained by slicing $N$ along
    this disc.  Since $M$ and $M'$ are homeomorphic to $N'$ and $N'$
    respectively, it follows that $M'$ is obtained by slicing
    along a properly embedded disc in $M$.
\end{proof}

\begin{figure}
\centering
\includegraphics{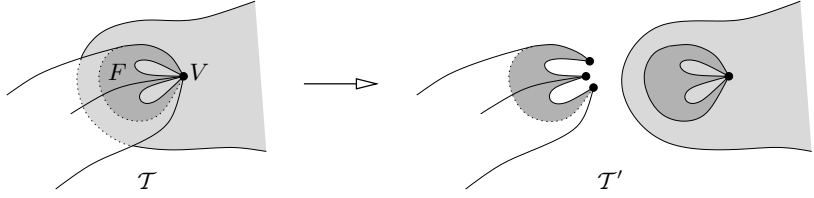}
\caption{Cutting along a face whose vertices are all identified}
\label{fig-prune-pinched}
\end{figure}

\begin{lemma} \label{l-prune-equator}
    Let $\tri$ be a triangulation of a bounded 3--manifold $M$.
    Let $e$ be a non-boundary edge in $\tri$ of degree one, let $\Delta$ be
    the (unique) tetrahedron containing $e$, and suppose that the edge
    opposite $e$ in $\Delta$ lies entirely within the boundary $\bdry M$,
    as illustrated in the leftmost diagram of Figure~\ref{fig-prune-equator}.

    \begin{figure}
    \centering
    \includegraphics{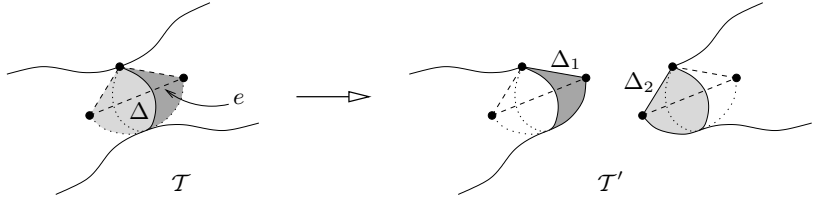}
    \caption{Cutting along a disc surrounding an edge of degree one}
    \label{fig-prune-equator}
    \end{figure}

    Suppose we replace $\Delta$ with two tetrahedra $\Delta_1$ and
    $\Delta_2$, each with two faces folded together to form an edge of
    degree one, as illustrated in the rightmost diagram of
    Figure~\ref{fig-prune-equator}.  Of the two portions of
    $\tri$ that were originally joined to $\Delta$ along the shaded faces,
    we join one of these portions to $\Delta_1$
    and the other to $\Delta_2$.  We leave the remaining faces of
    $\Delta_1$ and $\Delta_2$ as boundary faces.

    Then the result is a new triangulation $\tri'$ of some 3--manifold
    $M'$, which is homeomorphic to $M$ sliced along a properly
    embedded disc.
\end{lemma}

\begin{proof}
    In contrast to the previous result, this lemma contains no unusual cases.
    The edge opposite $e$ in $\Delta$ always bounds a properly embedded disc
    in $M$ (running directly through the centre of the tetrahedron $\Delta$),
    and so $\tri'$ triangulates a 3--manifold $M'$ that is obtained by
    slicing along this disc.
\end{proof}

In order to take full advantage of Lemmata~\ref{l-prune-face}
and~\ref{l-prune-equator}, it helps to have as many faces and edges of a
triangulation exposed to the boundary as possible.  The following operation
assists us in this regard.

\begin{lemma} \label{l-open-book}
    Let $\tri$ be a triangulation of a bounded 3--manifold $M$.
    Let $F$ be a non-boundary face of $\tri$, and supposed that precisely
    two of the three edges of $F$ lie within the boundary $\bdry M$,
    as illustrated in the leftmost diagram of Figure~\ref{fig-open-book}.
    Once again, neither the edges nor the vertices of $F$ are required to be
    distinct.

    \begin{figure}
    \centering
    \includegraphics{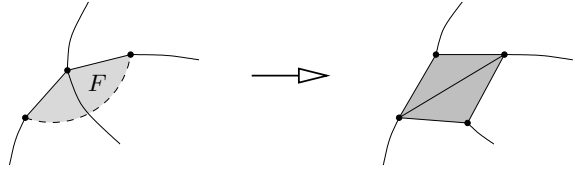}
    \caption{Performing a book opening move}
    \label{fig-open-book}
    \end{figure}

    Suppose that we unglue the two tetrahedra on either side of $F$,
    exposing these tetrahedra to the boundary as illustrated in
    the rightmost diagram of Figure~\ref{fig-open-book}.
    Then the result is a new triangulation $\tri'$ of the same 3--manifold $M$.

    We refer to this operation as a \emph{book opening move}.
\end{lemma}

\begin{proof}
    The proof is almost identical to Lemma~\ref{l-prune-face}, and we do
    not repeat the details.  The only differences are:
    \begin{itemize}
        \item Instead of slicing the manifold $M$ along a properly embedded
        disc, we slice it along a ``half-properly embedded'' disc.
        By this, we mean an embedded disc whose boundary consists of
        (i)~an arc in $\bdry M$, and (ii)~an arc in $\inter(M)$.
        Slicing along such a disc will never change the underlying
        manifold $M$.

        \item When we examine the neighbourhood of a vertex $V$,
        the intersection $F \cap \nbd(V)$ can include new types of
        discs that are ``half-properly embedded'' in $\nbd(V)$, as
        illustrated in Figure~\ref{fig-vertex-slice-half}.  Such discs are
        harmless however, and do not change the key fact that the resulting
        vertices in $\tri'$ all have 3--ball neighbourhoods.

        \begin{figure}
        \centering
        \includegraphics{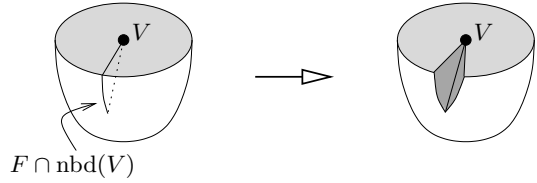}
        \caption{New types of discs in the intersection $F \cap \nbd(V)$}
        \label{fig-vertex-slice-half}
        \end{figure}
    \end{itemize}
    We refer to the proof of Lemma~\ref{l-prune-face} for the full details.
\end{proof}

We can now pull together all of these operations
to build a heuristic algorithm for detecting compressible surfaces.

\begin{algorithm}[Heuristic Pruning] \label{a-heuristic}
    Let $\tri$ be a triangulation of the closed 3--manifold $M$, and let
    $S$ be an embedded surface in $\tri$.  We can potentially show that
    $S$ is a compressible surface through the following procedure:
    \begin{enumerate}
        \item \label{en-prune-cutalong}
        Cut the triangulation $\tri$ along the surface $S$, and retriangulate
        the resulting bounded 3--manifold $M'$ (which may be disconnected).
        Let $\tri'$ denote this new bounded triangulation, let
        $\beta$ denote the number of boundary components,
        and let $\sigma$ denote the number of boundary spheres (so
        $\sigma \leq \beta$).

        \item \label{en-prune-simplify}
        Perform local simplification moves on $\tri'$ to reduce the
        number of tetrahedra (such as Pachner moves \cite{pachner91-moves}
        or boundary shellings).  There is no need to produce a minimal
        triangulation; it suffices to reach a point where there are no
        further immediate simplifications that can be done.

        \item \label{en-prune-openbook}
        Perform book opening moves on $\tri'$
        as described by Lemma~\ref{l-open-book} until no more can be done.

        \item \label{en-prune-search}
        Search for all locations within $\tri'$ at which the
        preconditions of Lemmata~\ref{l-prune-face} and~\ref{l-prune-equator}
        are satisfied.  That is, search for internal faces whose
        edges are all boundary, and search for internal degree one edges
        whose opposite edges are boundary.

        \item \label{en-prune-count}
        For each such location, temporarily perform the corresponding
        operation upon $\tri'$ that slices $M'$ along a properly
        embedded disc.
    \end{enumerate}
    If any of the sliced triangulations obtained in step~\ref{en-prune-count}
    still has $\beta$ boundary components \underline{or} still has $\sigma$
    boundary spheres, then the original surface $S$ is compressible
    in $M$.
\end{algorithm}

\begin{proof}
    Given Lemmata~\ref{l-prune-face}--\ref{l-open-book},
    the only part of this algorithm that remains to be proven is the final
    claim that, if some sliced triangulation from step~\ref{en-prune-count}
    has either $\beta$ boundary components or $\sigma$
    boundary spheres, then $S$ is compressible in $M$.

    Let $D$ be the properly embedded disc in $M'$ that we slice along in
    step~\ref{en-prune-count}.  If $D$ is \emph{not} a compressing disc
    for $S$ then $\bdry D$ bounds a disc in $\bdry M'$, whereupon
    slicing along $D'$ produces a new 2--sphere boundary component
    but otherwise leaves all existing boundary components unchanged.
    That is, we obtain $\beta+1$ boundary components, $\sigma+1$ of
    which are 2--spheres.

    Note that a compressing disc can produce a new boundary sphere
    without a new boundary component
    (for instance, slicing along the meridional disc of a solid torus),
    or a new boundary component without a new boundary sphere
    (for instance, slicing along a separating disc with non-trivial
    topology on each side).
    This is why we must count both boundary spheres \emph{and} boundary
    components in order to detect compressing discs.
\end{proof}

We finish this section with some notes regarding both the structure and
implementation of Algorithm~\ref{a-heuristic}.

\begin{itemize}
    \item The reason for steps~\ref{en-prune-simplify}
    and~\ref{en-prune-openbook} is to increase our chances of meeting
    the preconditions of Lemmata~\ref{l-prune-face} and~\ref{l-prune-equator}.
    In particular, simplifying the triangulation increases our chances
    of finding a region of $\tri'$ that is only ``one tetrahedron thick'',
    and book opening moves help expose more edges and vertices to the boundary.

    \item The precise local simplification moves of
    step~\ref{en-prune-simplify} are left up to the reader.  Many moves
    of this type are documented in the literature (particularly by
    authors involved in census enumeration); see
    \cite{burton04-facegraphs,matveev98-recognition}
    for some examples.

    The moves that we use in the following section with the {\swspace}
    include Pachner moves (also called bistellar moves
    \cite{pachner91-moves}), collapsing edges between distinct vertices,
    removing tetrahedra through boundary shellings, simplifying
    triangulations in the vicinity of low-degree edges and vertices,
    and the book opening move and its inverse (the \emph{book closing} move).

    \item Step~\ref{en-prune-cutalong}, in which we cut along the surface $S$
    in the triangulation $\tri$, causes a number of difficulties.
    The most severe problem is that it can generate a very large number
    of tetrahedra---in the case of the {\swspace}, we frequently find
    tetrahedra numbering in the thousands.
    It is therefore critical to have a simplification
    procedure that is both fast and effective.

    Moreover, cutting along a surface is messy for a programmer to implement,
    since tetrahedra can be subdivided into many different pieces of up
    to 11 distinct shapes (see Figure~\ref{fig-cut-shapes} for some
    examples).  Each of these shapes must be individually retriangulated
    (typically by the programmer as she implements the routine), and
    the code must then be able to automatically adjust the triangulations
    of these pieces so that the quadrilaterals, pentagons and hexagons
    on their boundaries can be glued together.

    \begin{figure}
    \centering
    \includegraphics[scale=0.7]{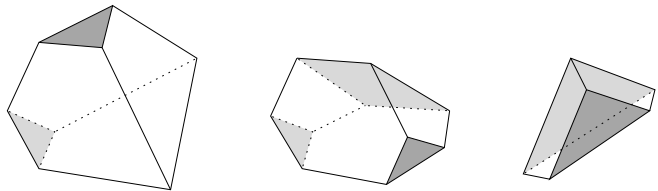}
    \caption{Sample pieces obtained by cutting along a normal surface}
    \label{fig-cut-shapes}
    \end{figure}

    We can make this cutting operation simpler if we first
    remove a small neighbourhood of each vertex of the original triangulation.
    This reduces the number of different shapes from 11 to 4,
    which is significantly easier for a programmer to manage.
    The boundaries of these pieces are also simpler to handle,
    with only quadrilaterals and hexagons to worry about.
    Of course we must not forget to glue the missing 3--balls back onto
    the boundary of the new triangulation once we are finished.
\end{itemize}

\section{The Weber-Seifert dodecahedral space} \label{s-swproof}

To conclude this paper, we apply the new Test~\ref{test-nonhaken}
to resolve an outstanding conjecture of Thurston.
The \emph{\swspace} is formed by identifying opposite faces of a
dodecahedron with a $3/10$ twist, and was one of the first known
examples of a hyperbolic 3--manifold \cite{weber33-dodecahedral}.
It was conjectured by Thurston that the {\swspace} is non-Haken
\cite{birman80-problems}, and here we prove this to be true.

By building on Haken's earlier work,
Jaco and Oertel gave an algorithm in 1984 to determine whether a
given 3--manifold $M$ has the Haken property \cite{jaco84-haken}.
This algorithm has since been improved by other authors
\cite{jaco02-algorithms-essential,tollefson98-quadspace},
though the basic framework remains the same:
\begin{enumerate}
    \item \label{en-jo-surfaces}
    Enumerate all vertex surfaces in some triangulation $\tri$ of $M$.
    \item \label{en-jo-discs}
    For each vertex surface $S \subset M$, test whether $S$ is incompressible:
    \begin{enumerate}[(i)]
        \item Cut $M$ along the surface $S$ and retriangulate.
        \item For each component $\tri'$ of the resulting triangulation,
        enumerate all fundamental normal surfaces in $\tri'$ and test
        whether any of these is a compressing disc.
    \end{enumerate}
\end{enumerate}
As noted in the introduction, all normal surfaces considered here are
embedded.  In particular, we do not consider the more general
case of immersed and/or singular surfaces within $\tri$.

In theory, proving the {\swspace} to be non-Haken should be a simple
matter of running the Jaco-Oertel algorithm.  However, this algorithm is
extremely slow in practice.  If $t$ is the number of tetrahedra in
$\tri$, then step~\ref{en-jo-surfaces} can grow exponentially slow in
$t$ and produce exponentially many surfaces, and each triangulation $\tri'$
in step~\ref{en-jo-discs} can contain exponentially many tetrahedra.
Even worse, each enumeration in step~\ref{en-jo-discs}(ii)
can grow exponentially slow in the size of $\tri'$, which becomes
\emph{doubly} exponential in $t$.  For these reasons, the Jaco-Oertel
algorithm has to date never been successfully applied to the {\swspace}.

Our proof begins in the same manner as step~\ref{en-jo-surfaces}
of the Jaco-Oertel algorithm---we triangulate the {\swspace}, and
then obtain a list of all vertex surfaces with the help of
recent developments in normal surface enumeration algorithms
\cite{burton09-convert,burton10-dd}.
However, the doubly exponential enumeration of step~\ref{en-jo-discs}
remains out of our reach, and so instead we use Theorem \ref{t-univ-disjoint}
and heuristic pruning to show that each of the surfaces in our list
has a compressing disc.

All of the computation in this proof was carried out using the
open-source software package {\regina} \cite{regina,burton04-regina}.
The supporting data for this proof (including the relevant triangulations
and normal surfaces) is contained in the
file {\swurl}, which readers can download from
the {\regina} website at \url{http://regina.sourceforge.net/data.html}.
%

\begin{definition}
    Let $\trisw$ denote the $23$--tetrahedron triangulation described in
    Table~\ref{tab-sw-triangulation}, which lists the pairwise
    identifications between the $4 \times 23$ faces of $23$ individual
    tetrahedra.  Here the tetrahedra are labelled $A,\ldots,W$
    and the vertices of each tetrahedron are numbered $1,\ldots,4$.

    Alternatively, we can describe $\trisw$ using the dehydration
    notation of Callahan, Hildebrand and Weeks.
    The dehydration string of $\trisw$ is
    \[\mathtt{xppphocgaeaaahimmnkontspmuuqrsvuwtvwwxwjjsvvcxxjjqattdwworrko},\]
    from which we can recover the full structure of $\trisw$ using the
    the rehydration procedure described in \cite{callahan99-cuspedcensus}.
\end{definition}

To read Table~\ref{tab-sw-triangulation}, each row gives the face
identifications for a single tetrahedron, and each column indicates one
of the four faces.  For instance, the cell in the bottom left corner indicates
that face 123 of tetrahedron $L$ is identified with face 214 of tetrahedron
$U$ (with vertices 1, 2 and 3 of $L$ identified with vertices 2, 1 and 4
of $U$ respectively).

\begin{table}
    \caption{The pairwise identifications of tetrahedron faces in the
        triangulation $\trisw$}
    \label{tab-sw-triangulation}
    \footnotesize \[ \begin{array}{c|r|r|r|r}
    & \mbox{123} & \mbox{124} & \mbox{134} & \mbox{234} \\
    \hline
    A & E: 123 & D: 124 & C: 134 & B: 234 \\
    B & H: 123 & G: 124 & F: 134 & A: 234 \\
    C & K: 123 & J: 124 & A: 134 & I: 234 \\
    D & M: 123 & A: 124 & H: 134 & L: 234 \\
    E & A: 123 & P: 124 & O: 134 & N: 234 \\
    F & M: 143 & Q: 124 & B: 134 & I: 134 \\
    G & R: 123 & B: 124 & N: 314 & M: 432 \\
    H & B: 123 & T: 124 & D: 134 & S: 234 \\
    I & O: 132 & K: 134 & F: 234 & C: 234 \\
    J & S: 123 & C: 124 & T: 134 & N: 231 \\
    K & C: 123 & V: 124 & I: 124 & U: 234 \\
    L & U: 214 & M: 142 & P: 132 & D: 234
    \end{array}
    \quad
    \begin{array}{c|r|r|r|r}
    & \mbox{123} & \mbox{124} & \mbox{134} & \mbox{234} \\
    \hline
    M & D: 123 & L: 142 & F: 132 & G: 432 \\
    N & J: 423 & U: 431 & G: 314 & E: 234 \\
    O & I: 132 & R: 234 & E: 134 & Q: 314 \\
    P & L: 143 & E: 124 & W: 134 & S: 431 \\
    Q & U: 213 & F: 124 & O: 324 & V: 134 \\
    R & G: 123 & T: 123 & W: 423 & O: 124 \\
    S & J: 123 & V: 123 & P: 432 & H: 234 \\
    T & R: 124 & H: 124 & J: 134 & W: 142 \\
    U & Q: 213 & L: 213 & N: 421 & K: 234 \\
    V & S: 124 & K: 124 & Q: 234 & W: 123 \\
    W & V: 234 & T: 243 & P: 134 & R: 341 \\
    \multicolumn{1}{c}{~} 
    \end{array} \]
\end{table}

\begin{lemma} \label{l-trisw}
    $\trisw$ is a triangulation of the {\swspace}.
\end{lemma}

\begin{proof}
    The triangulation $\trisw$ can be constructed as follows:
    \begin{enumerate}
        \item Build a regular dodecahedron by joining together twelve
        pentagonal cones, with the twelve apexes meeting at the centre of
        the dodecahedron and the twelve pentagonal bases forming the
        boundary of the dodecahedron.
        \item Triangulate each pentagonal cone with five tetrahedra, as
        illustrated in Figure~\ref{fig-sw-pentcone}.

        \begin{figure}
        \centering
        \includegraphics[scale=0.7]{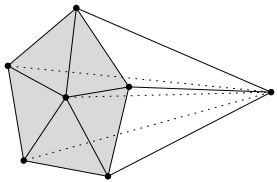}
        \caption{A pentagonal cone triangulated with five tetrahedra}
        \label{fig-sw-pentcone}
        \end{figure}

        \item Identify opposite faces of the dodecahedron with a $3/10$
        twist, giving a closed 3--manifold triangulation with 60 tetrahedra.
        \item \label{en-tri-simplify}
        Simplify this triangulation by first collapsing edges between
        distinct vertices, and then applying \mbox{3--2} Pachner moves (also
        called bistellar moves \cite{pachner91-moves}).  These operations are
        illustrated in Figure~\ref{fig-sw-simplify}.

        \begin{figure}
        \centering
        \includegraphics{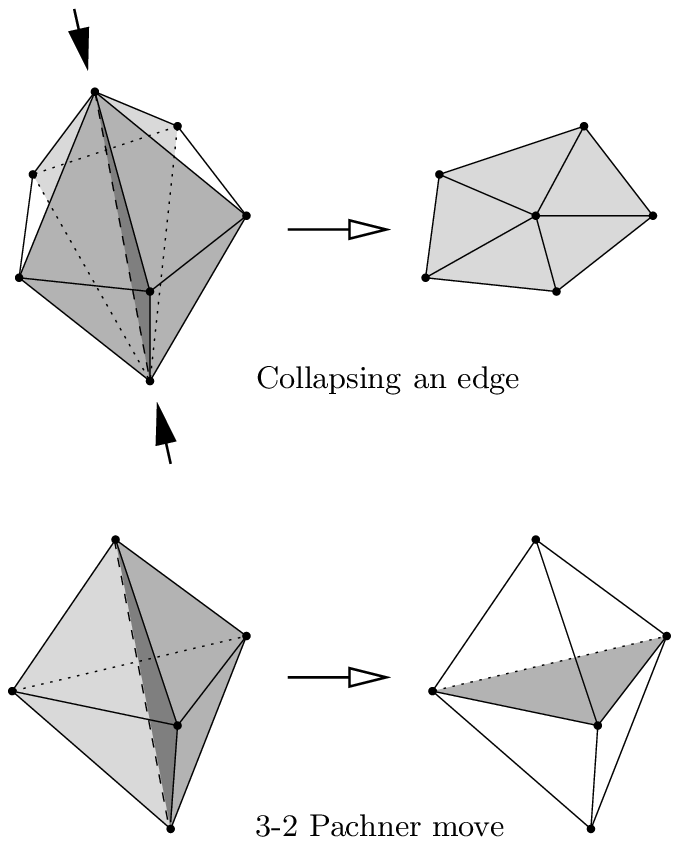}
        \caption{Moves to simplify the triangulation $\trisw$}
        \label{fig-sw-simplify}
        \end{figure}
    \end{enumerate}
    It is clear from this construction that $\trisw$ triangulates the
    {\swspace} as claimed.
\end{proof}

In fact, we conjecture that $\trisw$ is a \emph{minimal} triangulation
of the {\swspace}, i.e., that the space cannot be triangulated with
22~tetrahedra or fewer.  It should be noted that $\trisw$ is not the
only 23--tetrahedron triangulation of the {\swspace}; through a repeated
application of 2--3 and 3--2 Pachner moves we can obtain at least two
distinct\footnote{By \emph{distinct}, we mean that one triangulation
cannot be obtained from another simply by relabelling tetrahedra and
their vertices.}
alternatives, with the following dehydration strings:
\begin{equation}
\begin{split}
\mathtt{xppphjgbgaaaaginnlsnmupurtqsqurwwvvwvmwfcvuvkekaldmphextbvwfw} \\
\mathtt{xppfpnkalaaaamkionrmtnpoqrsutsqwuvvwwxdwvgtvqkpwxpxalnjcrkfns}
\end{split}
\label{eqn-sw-alt}
\end{equation}

Now that we are equipped with a triangulation, we embark on the first
step of the Jaco-Oertel algorithm---the enumeration of vertex
surfaces.  At the core of this enumeration is a linear programming
problem that takes place in a vector space of dimension $7t=161$
(where $t=23$ is the number of tetrahedra in $\trisw$).

A direct enumeration in $\R^{161}$ remains out of our reach
computationally, and so we take an indirect approach instead.
Tollefson \cite{tollefson98-quadspace} describes a smaller vector space
of dimension $3t$, in which we consider only the quadrilateral discs in
each normal surface.  We refer to this smaller vector space $\R^{3t}$ as
\emph{quadrilateral coordinates}, and in contrast we refer to the
original $\R^{7t}$ as used by Haken and then Jaco and Oertel as
\emph{standard coordinates}.

Our plan is (i)~to enumerate all vertex surfaces in quadrilateral
coordinates, and then (ii)~to convert this result into a list of all
vertex surfaces in standard coordinates.  Note that the latter step
is not just a matter of changing between coordinate systems, since the
``vertex surface'' property is not preserved between coordinate
systems---instead we must apply the complex (though extremely fast)
conversion procedure described in \cite{burton09-convert}.

In fact, step~(ii) is not necessary for the Jaco-Oertel algorithm, since
Tollefson proves that some $2$--sided incompressible surface must appear
as a vertex surface in quadrilateral coordinates, if such a
surface exists at all.  However, because our proof relies on
Theorem \ref{t-univ-disjoint}, we must work in
standard coordinates.  For us then, quadrilateral coordinates are
simply a means to an end.

By running the streamlined normal surface enumeration algorithm
described in \cite{burton10-dd}, we obtain the following result through
direct computation (recalling that vertex surfaces are defined here
to be $2$--sided and connected):

\begin{lemma} \label{l-enum-quad}
    The triangulation $\trisw$ has $698$ vertex surfaces in
    quadrilateral coordinates.  The genera of these surfaces are
    distributed according to the first row of Table~\ref{tab-surfacecounts}.
\end{lemma}

\begin{table}
\caption{Counting vertex surfaces according to genus}
\label{tab-surfacecounts}
\centering
\begin{tabular}{l|r|r|r|r|r|r|r|r|rr}
Genus & \phantom{00}0 & \phantom{00}1 & 2 & 3 & 4 & 5 &
    \phantom{00}6 & 7 & \phantom{00}8 \\
\cline{1-10}
Quad.~vertex surfaces & & 24 & 43 & 82 & 135 & 30 & & 300 & 36 \\
Standard vertex surfaces & 1 & 24 & 187 & 465 & 387 & 115 & 32 & 318 & 54 \\
\multicolumn{10}{c}{~} \\
Genus (ctd.) & 9 & 10 & 11 & 12 & 13 & 14 & 15 & 16 &
    \multicolumn{2}{|r}{Total} \\
\hline
Quad.~vertex surfaces & 36 & 12 & & & & & & &
    \multicolumn{2}{|r}{\textbf{698}} \\
Standard vertex surfaces & 54 & 36 & 30 & & 18 & 12 & 12 & 6 &
    \multicolumn{2}{|r}{\textbf{1751}} \\
\end{tabular}
\end{table}

The computation required to prove Lemma~\ref{l-enum-quad} is not
trivial---on a 2.4\,GHz Intel Core~2 CPU, the enumeration takes a little under
$5\frac12$ hours.  Without recent improvements to the normal surface
enumeration algorithm \cite{burton10-dd}, this computation could take
orders of magnitude longer, and without Tollefson's quadrilateral
coordinates it would remain completely infeasible.

Although Table~\ref{tab-surfacecounts} only lists the genus of each
surface, complete descriptions of all 698 surfaces can be found
in the file {\swurl}, as noted at the beginning of this section.

We now make our move into standard coordinates.  By running the 698 surfaces
of Lemma~\ref{l-enum-quad} through the quadrilateral-to-standard conversion
algorithm described in \cite{burton09-convert}, we obtain the following
result:

\begin{lemma} \label{l-enum-std}
    The triangulation $\trisw$ has $1751$ vertex surfaces in standard
    coordinates.  The genera of these surfaces are distributed according
    to the second row of Table~\ref{tab-surfacecounts}.
\end{lemma}

In contrast to the full enumeration in quadrilateral coordinates, the
conversion algorithm of \cite{burton09-convert} is extremely
fast, taking just over 1~second on the same 2.4\,GHz Intel Core~2 CPU.
As before, complete descriptions of all 1751 surfaces can be downloaded
in the file {\swurl}.

We pause here to make some observations about the low-genus surfaces in
our list.  By examining the individual normal discs that make up the spheres
and tori in Table~\ref{tab-surfacecounts}, we obtain the following result:

\begin{lemma} \label{l-links}
    The only vertex normal sphere in $\trisw$ is the frontier of a small
    regular neighbourhood of the single vertex of $\trisw$.
    Likewise, the only vertex normal tori in $\trisw$ are the frontiers
    of small regular neighbourhoods of the 24 edges of $\trisw$.
\end{lemma}

Using the nomenclature of Jaco and Rubinstein \cite{jaco03-0-efficiency},
these surfaces are called \emph{vertex links} and \emph{thin edge links}
respectively, and the triangulation $\trisw$ is both \emph{0--efficient}
and \emph{1--efficient} as a result.

Now that we have a full list of vertex surfaces at our disposal,
we can bring in the techniques of Theorem \ref{t-univ-disjoint} and
heuristic pruning to prove our final result.

\setcounter{c-save-section}{\arabic{section}}
\setcounter{c-save-theorem}{\arabic{theorem}}
\setcounter{section}{\arabic{c-sw-section}}
\setcounter{theorem}{\arabic{c-sw-theorem}}
\begin{theorem}
    The {\swspace} is non-Haken.
\end{theorem}
\setcounter{section}{\arabic{c-save-section}}
\setcounter{theorem}{\arabic{c-save-theorem}}

\begin{proof}
    Suppose the {\swspace} does contain a $2$--sided incompressible surface.
    Cases~2a and~2b of Theorem~\ref{t-univ-disjoint} are easily
    eliminated through a homology computation and Lemma~\ref{l-links},
    and so it follows from
    Theorem~\ref{t-univ-disjoint} that there must be two distinct,
    compatible, incompressible vertex surfaces $S_1,S_2$ in $\trisw$.
    We therefore run through our list of
    1751 vertex surfaces in search of such a pair $S_1,S_2$.

    We can eliminate the vertex linking sphere and the 24 vertex linking
    tori immediately.  Running Algorithm~\ref{a-heuristic} over the 1726
    remaining surfaces shows that 1710 of these contain a compressing disc.
    That is, heuristic pruning eliminates \emph{all but 16} of these
    vertex surfaces.  Those surfaces that remain are
    summarised in Table~\ref{tab-surfacecounts-pruned} (once again,
    see {\swurl} for their full descriptions).

    \begin{table}
    \caption{The 16 vertex surfaces that remain after heuristic pruning}
    \label{tab-surfacecounts-pruned}
    \centering
    \begin{tabular}{l|r|r|r|r|r|r}
    Genus & 4 & 5 & 7 & 8 & 9 & Total \\
    \hline
    Surfaces remaining & 1 & 3 & 8 & 1 & 3 & \textbf{16}
    \end{tabular}
    \end{table}

    It follows that, if they exist at all, the surfaces $S_1$ and $S_2$ must
    belong to this smaller list.  However, comparing quadrilateral types for
    all $\binom{16}{2}$ pairs shows that no two of these surfaces
    are compatible, and so by Theorem~\ref{t-univ-disjoint} the
    {\swspace} cannot be Haken.
\end{proof}

We finish with a handful of observations regarding the different
elements used in the proof of Theorem~\ref{t-sw}.

\begin{itemize}
    \item
    It is mentioned earlier that the number of vertex surfaces can
    grow exponentially in the number of tetrahedra $t$.  The best
    theoretical bounds known to date are based on
    the upper bound theorem of McMullen \cite{mcmullen70-ubt}, yielding
    theoretical limits of $O(4^t)$ in quadrilateral coordinates and
    $O(15^t)$ in standard coordinates \cite{burton10-dd}.
    It is therefore surprising in our case with $t=23$ to find just
    698 and 1751 vertex surfaces respectively.  Such
    discrepancies between theory and practice are common, and are
    discussed in greater detail in the paper \cite{burton10-complexity}.

    \item
    We can recall from Section~\ref{s-heuristic} that heuristic pruning
    involves two distinct tests: one for internal faces
    with three boundary edges (Lemma~\ref{l-prune-face}) and
    one for discs surrounding edges of degree one
    (Lemma~\ref{l-prune-equator}).  It is worth comparing the relative
    effectiveness of these tests.

    Of the 1726 vertex surfaces upon which we attempt heuristic
    pruning, 1695 can be eliminated using Lemma~\ref{l-prune-face}
    but only 88 can be eliminated through Lemma~\ref{l-prune-equator}.
    These are success rates of approximately $98\%$ and $5\%$
    respectively.  It appears therefore that testing for faces with
    three boundary edges is significantly more powerful in practice.

    \item
    Tollefson proves that quadrilateral coordinates are
    sufficient for running the original Jaco-Oertel algorithm
    \cite{tollefson98-quadspace}, whereas in this paper we use standard
    coordinates instead.  It is worth noting that this choice does not
    lead to any significant loss of efficiency or power:
    \begin{itemize}
        \item Assuming that we already have a list of
        vertex surfaces in quadrilateral coordinates, creating a list
        of vertex surfaces using the conversion algorithm of
        \cite{burton09-convert} is extremely fast, taking only a matter of
        seconds of processing time.

        \item Although we have more surfaces to deal with in standard
        coordinates (1751 instead of 698), heuristic pruning eliminates
        these differences entirely.  That is, applying
        heuristic pruning to the vertex surfaces in \emph{quadrilateral}
        coordinates leaves us with precisely the same 16 surfaces that
        we describe in Table~\ref{tab-surfacecounts-pruned}.
        Similar behaviour is seen when working with the alternate
        triangulations described by (\ref{eqn-sw-alt}).
    \end{itemize}

    \item
    Although the triangulation $\trisw$ was chosen arbitrarily, in
    hindsight this was a fortuitous choice.  If we attempt to apply the
    method used in Theorem~\ref{t-sw} to either of the alternative
    triangulations described by (\ref{eqn-sw-alt}), we do not arrive at
    a conclusive proof.

    Specifically, if we (i)~eliminate vertex and thin edge links,
    (ii)~eliminate surfaces through heuristic pruning, and then
    (iii)~eliminate surfaces without compatible partners according to
    Theorem~\ref{t-univ-disjoint}, some surfaces still remain.
    For the first alternative we are left with one compatible pair of
    genus~7 surfaces, and for the second alternative we are left with
    three compatible pairs of genus~7 surfaces.  All of these leftover
    surfaces are vertex surfaces in both standard and quadrilateral
    coordinates.  They
    can eventually be eliminated, but only with additional manipulation
    of the corresponding bounded triangulations.

    \item
    It is interesting to compare the relative power and efficiency of
    Theorem \ref{t-univ-disjoint} and heuristic pruning as individual
    techniques.  Ignoring the vertex link and thin edge links, consider the
    remaining 1726 vertex surfaces in $\trisw$.
    As we have seen already, heuristic pruning
    alone eliminates 1710 of these 1726 surfaces (around $99\%$).
    On the other hand, if we use Theorem \ref{t-univ-disjoint} as a
    filtering tool (by removing all surfaces with no compatible partner
    as required by Theorem~\ref{t-univ-disjoint}), we can eliminate 1227
    of these 1726 surfaces (around $71\%$).

    Although this suggests that Theorem \ref{t-univ-disjoint} is less
    effective as a filtering tool than heuristic pruning, it is
    significantly faster to use---filtering by Theorem \ref{t-univ-disjoint}
    takes just a few seconds, whereas running all 1726 surfaces
    through heuristic pruning takes about 40~minutes (primarily because
    we must cut along each surface, which can produce triangulations
    with thousands of tetrahedra to simplify and test).  It follows that
    Theorem \ref{t-univ-disjoint} could be used as a very fast initial filter,
    leaving a smaller set of surfaces to run through the more expensive
    heuristic pruning.

    A more sophisticated variant of this idea is to repeatedly call upon
    Theorem~\ref{t-univ-disjoint} \emph{throughout} the heuristic
    pruning process.  That is, every time a surface $S$ is eliminated
    through heuristic pruning, we immediately eliminate every other
    surface that has $S$ as its only compatible partner.  Although this
    should further improve the efficiency of elimination, the final result
    (i.e., the set of leftover surfaces) will of course remain the same.
\end{itemize}

Finally, we note that it is possible to prove the {\swspace} to be
non-Haken without employing Theorem \ref{t-univ-disjoint} at all,
although the relevant computations
require significantly more human intervention.  If we begin with the
first alternative triangulation of (\ref{eqn-sw-alt}), we obtain
1909 vertex surfaces.  With heuristic pruning this list reduces
to just nine surfaces: eight surfaces $S_1,\ldots,S_8$ of genus~7,
and one surface $S_9$ of genus~8.

For each genus~7 surface $S_i$ ($1 \leq i \leq 8$), cutting along $S_i$
gives at least one bounded triangulation $\tri_i$ with free fundamental group.
By repeatedly applying Pachner moves, we can recognise the underlying
manifold as a genus~7 handlebody, showing the original
surface $S_i$ to be compressible.

The final surface $S_9$ is more difficult to deal with.  We cut along
$S_9$ to obtain bounded triangulations $\tri_9$ and $\tri_9'$; although
we are not able to identify either component, with sufficiently many
Pachner moves we can nevertheless manufacture a
compressing disc in the form described by Lemma~\ref{l-prune-face}.

%
%

\providecommand{\bysame}{\leavevmode\hbox to3em{\hrulefill}\thinspace}
\providecommand{\MR}{\relax\ifhmode\unskip\space\fi MR }
\providecommand{\MRhref}[2]{%
  \href{http://www.ams.org/mathscinet-getitem?mr=#1}{#2}
}
\providecommand{\href}[2]{#2}

\end{document}